\documentclass[12pt,reqno]{amsart}
\usepackage{indentfirst, amssymb,amsmath,amsthm}
\newtheorem{theo}{Theorem}[section]

\newtheorem{pro}{Proposition}[section]
\newtheorem{cor}{Corollary}[section]
\newtheorem{defi}{Definition}[section]

\newcommand{\be}{\begin{equation}}
\newcommand{\ee}{\end{equation}}
\newcommand{\beas}{\begin{eqnarray*}}
\newcommand{\eeas}{\end{eqnarray*}}
\newcommand{\bea}{\begin{eqnarray}}
\newcommand{\eea}{\end{eqnarray}}

\numberwithin{equation}{section}

\begin{document}

\setlength{\unitlength}{1mm} \baselineskip .45cm
\setcounter{page}{1}
\pagenumbering{arabic}
\title[]
{Conformal vector fields on almost Kenmotsu manifolds}

\author[]{Uday Chand De, Arpan Sardar and Krishnendu De$^{*}$}

\address{Uday Chand De \newline  Department of Pure Mathematics,\newline University of Calcutta,
\newline 35, Ballygunge Circular Road,\newline Kolkata 700019,\newline West Bengal, India.\newline ORCID iD: https://orcid.org/0000-0002-8990-4609}
 \email {uc$\_$de@yahoo.com}

\address {\newline Arpan Sardar \newline Department of Mathematics \newline University of Kalyani \newline Kalyani 741235
\newline Nadia, West Bengal, India}
\email {arpansardar51@gmail.com}
\address
 { Krishnendu De \newline Department of Mathematics,\newline Kabi Sukanta Mahavidyalaya,\newline
The University of Burdwan. \newline Bhadreswar, P.O.-Angus, Hooghly,\newline
Pin 712221, West Bengal, India.\newline ORCID iD: https://orcid.org/0000-0001-6520-4520}
\email{krishnendu.de@outlook.in }

\begin{abstract}
In this paper, first we consider that the conformal vector field $\mathbf{X}$ is identical with the Reeb vector field $\varsigma$ and next, assume that $\mathbf{X}$ is pointwise collinear with 
$\varsigma$, in both cases it is shown that the manifold $\mathbf{N}^{2m+1}$ becomes a Kenmotsu manifold and $\mathbf{N}^{2m+1}$ is locally a warped product $\mathbf{N}' \times_{f} \mathbf{M}^{2m}$, where $\mathbf{M}^{2m}$ is an almost K\"ahler manifold, $\mathbf{N}'$ is an open interval with coordinate t, and $f = ce^{t}$ for some positive constant c. Beside these, we prove that if a $(\verb"k",\boldsymbol{\mu})'$-almost Kenmotsu manifold admits a Killing vector field $\mathbf{X}$, then either it is locally a warped product of an almost K\"ahler manifold and an open interval or $\mathbf{X}$ is a strict infinitesimal contact transformation. Furthermore, we also investigate $\boldsymbol{\eta}$-Ricci-Yamabe soliton with conformal vector fields on  $(\verb"k",\boldsymbol{\mu})'$-almost Kenmotsu manifolds and finally, we construct an example.

\end{abstract}

\maketitle
\footnotetext{\\
\keywords{Key words and phrases : conformal vector fields, infinitesimal strict contact transformation, $\boldsymbol{\eta}$-Ricci-Yamabe solitons, almost Kenmotsu manifolds, $(k,\boldsymbol{\mu})'$-almost Kenmotsu manifolds.\\
\thanks{$^{*}$ Corresponding author}
}
}

\section{\textsf{Introduction}}

In a contact manifold, the existence of conformal vector fields (briefly, $CVF$) can be explained beautifully by some intrinsic properties of contact manifolds. For example, it is well circulated that \cite{yan} a Riemannian manifold $\mathbf{N}^m$ of dimension $m$ is conformally flat if it concedes a maximal, that is, $\frac{(m+1)(m+2)}{2}$-parameter group of conformal motions. A conformally flat Sasakian manifold is also known to be of constant curvature $1$ \cite{okm}. In \cite{okm1}, Okumura established that for $m>3$, a connected complete Sasakian manifold equipped with a conformal motion is isometric to sphere. After that, Sharma \cite{sharma1} extended his research on contact manifolds that admit conformal motions to $N(\verb"k")$-contact metric manifolds. Later, the idea of the holomorphic planar $CVF$ in Hermitian manifolds \cite{sharma2} was also introduced by Sharma. Planar $CVF$ have also been studied in \cite{sharma2}.\par

A vector field $\mathbf{X}$ on $\mathbf{N}^{2m+1}$ satisfying the equation
\begin{equation}\label{1.1}
\pounds_\mathbf{X} g = 2 \rho g,
\end{equation}
$\rho$ being a smooth function and $\pounds$ is the Lie-derivative, is called a  $CVF$. If $\mathbf{X}$ is not Killing, it is termed as non-trivial. If $\rho$ vanishes, then the $CVF$ $\mathbf{X}$ is named Killing. $\mathbf{X}$ is called homothetic, if $\rho$ is constant.  $CVF$ have been studied by many authors such as (\cite{de1}, \cite{des1}, \cite{des7}, \cite{oba3}, \cite{sharma1}-\cite{sharma3}) and many others.\par

If $r, R, S$ indicate the scalar curvature, the curvature tensor and the Ricci tensor, respectively, then the $CVF$ $\mathbf{X}$ satisfies the followings\cite{yan}
\begin{equation}\label{1.2}
(\pounds_\mathbf{X} \nabla)(U_1,V_1) = (U_1\rho)V_1 - (V_1\rho)U_1 -g(U_1,V_1)D\rho,
\end{equation}
\begin{eqnarray}\label{1.3}
(\pounds_\mathbf{X} R)(U_1,V_1)W_1 &=& g(\nabla_{U_1} D\rho,W_1)V_1 - g(\nabla_{V_1} D\rho,W_1)U_1\\ \nonumber
&& + g(U_1,W_1)\nabla_{V_1} D\rho - g(V_1,W_1)\nabla_{U_1} D\rho,
\end{eqnarray}
\begin{equation}\label{1.4}
(\pounds_\mathbf{X} S)(U_1,V_1) = -(2m-1)g(\nabla_{U_1} D\rho,V_1) - (\bigtriangleup \rho)g(U_1,V_1),
\end{equation}
\begin{equation}\label{1.5}
\pounds_\mathbf{X} r = -4m(\bigtriangleup\rho) -2r\rho
\end{equation}
for all vector fields $U_1, V_1, W_1$ on $\mathbf{N}^{2m+1}$, where $D{\rho}$ and $\bigtriangleup {\rho} = div D{\rho}$ respectively denote the gradient and Laplacian of $\rho$. \par

\begin{defi}
A vector field $\mathbf{X}$ satisfying the relation
\begin{equation}\label{1.6}
\pounds_\mathbf{X} \boldsymbol{\eta} = \sigma \boldsymbol{\eta},
\end{equation}
 $\sigma$ being a scalar function, is named an infinitesimal contact transformation. It is named as infinitesimal strict contact transformation, if $\sigma$ vanishes identically.
\end{defi}

Sharma and Vrancken\cite{sharma3} generalized the theorem of Tanno\cite{tanno} and proved:
\begin{theo}
The $CVF$ on $\mathbf{N}$ is an infinitesimal automorphism of $\mathbf{N}$ if it is an infinitesimal contact transformation.
\end{theo}

Now, for dimension $>$ 3, we recollect the subsequent outcome on a Sasakian manifold stated by Okumura\cite{okm1}.

\begin{theo}
If $\mathbf{N}^{2m+1},\, m> 1$ be a Sasakian manifold with a non-Killing  $CVF$ $\mathbf{X}$, then $\mathbf{X}$ is special concircular. Moreover, $\mathbf{N}^{2m+1}$ is isometric to $\mathbb{S}^{2m+1}$ (unit sphere) if it is connected and complete as well.
\end{theo}

Again, Sharma and Blair(\cite{sharma1}) proved the following Theorem to generalize the foregoing result.
\begin{theo}
On a $(\verb"k",0)$-contact manifold $\mathbf{N}^{2m+1}$, let $\mathbf{X}$ be a non-Killing $CVF$. $\mathbf{N}^{2m+1}$ is Sasakian and $\mathbf{X}$ is concircular for $m > 3$, hence if $\mathbf{N}^{2m+1}$ is complete, it is isometric to $\mathbb{S}^{2m+1}$. Also, $\mathbf{N}^{2m+1}$ is either Sasakian or flat for $m = 1$.
\end{theo}

To generalized the above result, Sharma and Vrancken\cite{sharma3} established the following theorem:
\begin{theo}
Let $\mathbf{N}^{2m+1}$ be a $(\verb"k",\boldsymbol{\mu})$-contact metric manifold which admits $\mathbf{X}$, a non-Killing $CVF$. For $m > 3$,\\
(i) $\mathbf{N}^{2m+1}$ is Sasakian and $\mathbf{X}$ is concircular, also $\mathbf{N}^{2m+1}$ is isometric to $\mathbb{S}^{2m+1}$ if it is complete, or\\
(ii) $\verb"k" =-m-1$ and $\boldsymbol{\mu}= 1$. Also, $\mathbf{N}^{2m+1}$ is isometric $S^{2m+1}$ if it is compact.

\end{theo}
Very recently, De, Suh and Chaubey\cite{de1} studied  $CVF$ on almost co-K\"ahler manifolds. In 2022 \cite{wana1}, Wang investigated almost Kenmotsu $(\verb"k",\boldsymbol{\mu})'$-manifolds with $CVF$ in dimension three.\par

Inspired by the foregoing studies we are interested to investigate  $CVF$ on almost Kenmotsu manifolds (briefly, akm).\par

The following is how the current article is structured:\par
Following the preliminaries, we investigate $CVF$ in akm in Section 3. In next Section, we investigate a $(\verb"k",\boldsymbol{\mu})'$-akm admitting a $CVF$ $\mathbf{X}$. Section 5, concerns with the study of $\boldsymbol{\eta}$-Ricci-Yamabe soliton with $CVF$ on  $(\verb"k",\boldsymbol{\mu})'$-akm. Finally, we build an example.
\vspace{.8cm}

\section{\textsf{Preliminaries}}

A Riemannian manifold $(\mathbf{N}^{2m+1},g)$ is caled an almost contact metric manifold \cite{blair} if it admit a tensor field $\varphi$ ($(1,1)$-type), a vector field $\varsigma$ and a 1-form $\boldsymbol{\eta}$ obeying
\begin{equation}\label{2.1}
\varphi ^2 U_1 = -U_1 + \boldsymbol{\eta}(U_1)\varsigma,\hspace{.3cm} \boldsymbol{\eta}(\varsigma) = 1,
\end{equation}
\begin{equation}\label{2.2}
g(\varphi U_1, \varphi V_1) = g(U_1,V_1) - \boldsymbol{\eta}(U_1)\boldsymbol{\eta}(V_1)
\end{equation}
for all vector fields $U_1,V_1$. The vector field $\varsigma$ is called the Reeb or characteristic vector field.\par
If the Nijenhuis tensor of $\varphi$ \cite{blair} disappears, an almost contact metric manifold is called normal.
If in an almost contact metric manifold, $d\boldsymbol{\eta} = 0$ and $d\varphi = 2\boldsymbol{\eta} \wedge \varphi$, $\varphi (U_1,V_1) = g(U_1,\varphi V_1)$, then it is akm. A normal almost akm is a Kenmotsu manifold. In an akm the relation
\begin{equation}\label{2.3}
(\nabla_{U_1} \varphi)V_1 = g(\varphi U_1, V_1)\varsigma - \boldsymbol{\eta}(V_1)\varphi U_1
\end{equation}
holds. Also the following formulas hold in akm (\cite{dil1},\cite{dil2})
\begin{equation}\label{2.4}
h\varphi + \varphi h = 0,\hspace{.3cm}h\varsigma = h'\varsigma = 0,\hspace{.3cm}  tr(h) = tr(h') = 0,
\end{equation}
\begin{equation}\label{2.5}
\nabla_{U_1} \varsigma = U_1 - \boldsymbol{\eta}(U_1)\varsigma + h'U_1,
\end{equation}
where $h=\frac{1}{2}\pounds_\varsigma \varphi$ and $h' = h \circ \varphi$.

In an akm if the characteristic vector field $\varsigma$ belongs to $(\verb"k",\boldsymbol{\mu})'$-nullity distribution, that is,
\begin{equation}\label{2.6}
R(U_1,V_1)\varsigma = \verb"k"(\boldsymbol{\eta}(V_1)U_1-\boldsymbol{\eta}(U_1)V_1) + \boldsymbol{\mu}(\boldsymbol{\eta}(V_1)h'U_1 - \boldsymbol{\eta}(U_1)h'V_1),
\end{equation}
$\verb"k"$ and $\boldsymbol{\mu}$ are constants, then it is called a $(\verb"k",\boldsymbol{\mu})'$-akm\cite{dil1}. In a $(\verb"k",\boldsymbol{\mu})'$-akm \cite{dil1}:
\begin{equation}\label{2.7}
h'^2 U_1 = (\verb"k"+1)\varphi^2 U_1
\end{equation}
 and $\boldsymbol{\mu} = -2$. Equation (\ref{2.7}) reflects that $h'=0$  if and only if $\verb"k"+1=0$ whereas $h'\neq 0$ if and only if $\verb"k"+1 < 0$. The equation (\ref{2.6}) yields
\begin{equation}\label{2.8}
R(\varsigma,U_1)V_1 = \verb"k"(g(U_1,V_1)\varsigma - \boldsymbol{\eta}(V_1)U_1) - 2(g(h'U_1,V_1)\varsigma-\boldsymbol{\eta}(V_1)h'U_1)
\end{equation}
for any  $U_1,V_1 \in \chi(\mathbf{N})$.\\
\begin{pro}(\cite{wan})
In a $(\verb"k",\boldsymbol{\mu})'$-akm  the $(1,1)$-Ricci tensor $Q$ is given by
\begin{equation}\label{2.9}
QU_1 = -2mU_1 + 2m(\verb"k"+1)\boldsymbol{\eta}(U_1)\varsigma -2mh'U_1
\end{equation}
for $\verb"k"+1<0$, where $Q$ is defined by $S(U_1,V_1) = g(QU_1,V_1)$ and 
$r= 2m(\verb"k"-2m)$.
\end{pro}

\begin{pro}(\cite{dil2})
An akm is a Kenmotsu manifold if and only if h = 0.
\end{pro}

\begin{pro}(\cite{dil2})
 If $\mathbf{N}^{2m+1}$ be an akm with $h = 0$, then $\mathbf{N}^{2m+1}$ is locally a warped product $\mathbf{N}' \times_{f} \mathbf{M}^{2m}$, in which $\mathbf{M}^{2m}$ indicate an almost K\"ahler manifold, with coordinate $t$, $\mathbf{N}'$ being the open interval and $f = ce^{t}$ for some $c$ ( positive constant).
\end{pro}

 \begin{defi}
 An akm $\mathbf{N}^{2m+1}$ is named an $\boldsymbol{\eta}$-Einstein manifold if its Ricci tensor $S$ fulfills
 \begin{equation}\nonumber
 S(U_1,V_1) = a_1g(U_1,V_1) + b_1\boldsymbol{\eta}(U_1)\boldsymbol{\eta}(V_1),
  \end{equation}
 where $a_1,\, b_1$ are scalars of which $b_1\neq 0$.
\end{defi}
In \cite{ken}, the author proved that an $\boldsymbol{\eta}$-Einstein Kenmotsu manifold is an Einstein manifold, provided $b_1$ = constant (or, $a_1$ = constant). Also, Pastore and Saltarelli\cite{pas} proved that $\boldsymbol{\eta}$-Einstein $(\verb"k",\boldsymbol{\mu})$-akm is an Einstein manifold for any one of $a_1$ or $b_1$ is constant. Again, Mandal and De\cite{man1} studied the above result in $(\verb"k",\boldsymbol{\mu})'$-akm. They explain the following:
\begin{pro}
	An $\boldsymbol{\eta}$-Einstein $(\verb"k",\boldsymbol{\mu})'$-akm becomes an Einstein manifold, provided 
$a_1$ or $b_1$ is constant.
\end{pro}
\vspace{.6cm}

\section{\textsf{Conformal vector fields on almost Kenmotsu manifolds}}

Let us assume that the Reeb vector field $\varsigma$ be a  $CVF$ on $\mathbf{N}^{2m+1}$. Then equation (\ref{1.1}) implies
\begin{equation}\label{3.1}
(\pounds_{\varsigma} g)(U_1,V_1) = 2\rho g(U_1,V_1),
\end{equation}
which means that
\begin{equation}\label{3.2}
g(\nabla_{U_1} \varsigma,V_1) + g(U_1,\nabla_{V_1} \varsigma) = 2\rho g(U_1,V_1).
\end{equation}
Using (\ref{2.5}) in (\ref{3.2}), we obtain
\begin{equation}\label{3.3}
\rho g(U_1,V_1) = [g(U_1,V_1)-\boldsymbol{\eta}(U_1)\boldsymbol{\eta}(V_1) + g(h'U_1,V_1)].
\end{equation}
Setting $U_1=V_1=\varsigma$ in the foregoing equation entails that
\begin{equation}\label{3.5}
\rho = 0.
\end{equation}
Applying (\ref{3.5}) in (\ref{3.1}), we obtain
\begin{equation}
(\pounds_{\varsigma} g)(U_1,V_1) = 0.
\end{equation}
Hence, the foregoing equation implies that $\varsigma$ is a Killing. Thus, we acquire $h=0$. Therefore, from Proposition 2.2 and 2.3, we have:
\begin{theo}
If the Reeb vector field $\varsigma$ of $\mathbf{N}^{2m+1}$ is a $CVF$, then $\mathbf{N}^{2m+1}$ becomes a Kenmotsu manifold and $\mathbf{N}^{2m+1}$ is also locally a warped product $\mathbf{N}' \times_{f} \mathbf{M}^{2m}$, where $\mathbf{M}$, $\mathbf{N}'$ and $f$ are defined earlier.
\end{theo}

\vspace{.4cm}

Suppose $\mathbf{X}=b\varsigma$, where $b$ is a smooth function on $\mathbf{N}^{2m+1}$. Then, using (\ref{1.1}), we get
\begin{equation}\label{3.7}
(\pounds_{b\varsigma} g)(U_1,V_1) = 2\rho g(U_1,V_1),
\end{equation}
which implies
\begin{equation}\label{3.8}
g(\nabla_{U_1} b\varsigma, V_1) + g(U_1, \nabla_{V_1} b\varsigma) = 2\rho g(U_1,V_1).
\end{equation}
Using (\ref{2.5}) in the foregoing equation gives
\begin{equation}\label{3.9}
(U_1 b)\boldsymbol{\eta}(V_1) + (V_1 b)\boldsymbol{\eta}(U_1) + 2b[g(U_1,V_1)-\boldsymbol{\eta}(U_1)\boldsymbol{\eta}(V_1)+ g(h'U_1,V_1)] = 2\rho g(U_1,V_1).
\end{equation}
Contracting the above equation we obtain
\begin{equation}\label{3.10}
\varsigma b = \rho\,(2m+1)-2mb.
\end{equation}
Replacing $U_1$ by $\varphi U_1$ in (\ref{3.9}), we get
\begin{equation}\label{3.11}
(\varphi U_1 b)\boldsymbol{\eta}(V_1) + 2b\,[g(\varphi U_1,V_1)-g(h U_1,V_1)]=2\rho g(\varphi U_1, V_1).
\end{equation}
Putting $V_1 =\varsigma$ in (\ref{3.11}), we get
\begin{equation}\label{3.12}
(\varphi U_1) b =0.
\end{equation}
The above two equations implies
\begin{equation}\label{3.13a}
b\,[g(\varphi U_1,V_1)-g(h U_1,V_1)]=\rho g(\varphi U_1, V_1).
\end{equation}
Replacing $U_1$ by $V_1$ and $V_1$ by $U_1$ in (\ref{3.13a}), we get
\begin{equation}\label{3.14b}
b\,[-g(\varphi U_1,V_1)-g(h U_1,V_1)]=-\rho g(\varphi U_1, V_1).
\end{equation}
In view of the above two equations, we get
\begin{equation}\label{3.15c}
(b-\rho)g(\varphi U_1,V_1)=0,
\end{equation}
which implies
\begin{equation}\label{3.16d}
b=\rho.
\end{equation}
Using (\ref{3.16d}) in (\ref{3.14b}), we infer
\begin{equation}
b\, g(h U_1,V_1)=0,
\end{equation}
which implies either $b=0$ or, $b\neq 0$.\par
Case I: If $b=0$, then (\ref{3.16d}) implies $\rho =0$. Hence, $\mathbf{X}$ is Killing, which means that $\varsigma$ is Killing.\\
Case II: If $b\neq 0$, then $g(hU_1,V_1)=0$. Hence $h=0$, which implies $\varsigma$ is Killing.\par
Hence, both cases implies $\varsigma$ is Killing.  Hence using Proposition 2.2 and 2.3, we provide:
\begin{theo}
In $\mathbf{N}^{2m+1}$, if a $CVF$ $\mathbf{X}$  is pointwise collinear with 
$\varsigma$, then $\mathbf{N}^{2m+1}$ becomes a Kenmotsu manifold and $\mathbf{N}^{2m+1}$ is locally a warped product $\mathbf{N}' \times_{f} \mathbf{M}^{2m}$.
\end{theo}

\vspace{.6cm}

\section{\textsf{Conformal vector fields on $(k,\boldsymbol{\mu})'$-akm}}

Suppose the vector field $\mathbf{X}$ in $\mathbf{N}^{2m+1}$ is Killing. Then
\begin{equation}\label{3.13}
(\pounds_\mathbf{X} g)(U_1,V_1) =0
\end{equation}
and
\begin{equation}\label{3.14}
(\pounds_\mathbf{X} S)(U_1,V_1) = 0.
\end{equation}
Definition of Lie-derivative infers that
\begin{equation}\label{3.15}
(\pounds_\mathbf{X} \boldsymbol{\eta})U_1 = \pounds_\mathbf{X} \boldsymbol{\eta}(U_1) - \boldsymbol{\eta}(\pounds_\mathbf{X} U_1).
\end{equation}
Equation (\ref{3.13}) implies
\begin{equation}\label{3.16}
(\pounds_\mathbf{X} \boldsymbol{\eta})U_1 = g(\pounds_\mathbf{X} \varsigma,U_1).
\end{equation}
Also, we have
\begin{equation}\label{3.17}
\boldsymbol{\eta}(\pounds_\mathbf{X} \varsigma) =0 \hspace{.4cm} and \hspace{.4cm} (\pounds_\mathbf{X} \boldsymbol{\eta})\varsigma = 0.
\end{equation}
From (\ref{2.9}), we obtain
\begin{equation}\label{3.17.1}
S(U_1,V_1) = -2m g(U_1,V_1)+2m(\verb"k"+1)\boldsymbol{\eta}(U_1)\boldsymbol{\eta}(V_1)-2mg(h'U_1,V_1).
\end{equation}
Now, we take Lie-derivative of 
(\ref{3.17.1}) along 
$\mathbf{X}$ and get
\begin{eqnarray}\label{3.18}
(\pounds_\mathbf{X} S)(U_1,V_1) &=&-2m(\pounds_{\mathbf{X}}g)(U_1,V_1) \\ \nonumber
&&+2m(\verb"k"+1)[((\pounds_{\mathbf{X}}\boldsymbol{\eta})U_1) \boldsymbol{\eta}(V_1)+((\pounds_{\mathbf{X}}\boldsymbol{\eta})V_1) \boldsymbol{\eta}(U_1)]\\ \nonumber
&&-2m[(\pounds_{\mathbf{X}}g)(h'U_1,V_1)+g((\pounds_{\mathbf{X}}h')U_1,V_1)].
\end{eqnarray}
Using (\ref{3.13}) and (\ref{3.14}) in (\ref{3.18}), we provide
\begin{eqnarray}\label{3.19}
&&(\verb"k"+1)[((\pounds_{\mathbf{X}}\boldsymbol{\eta})U_1) \boldsymbol{\eta}(V_1)+((\pounds_{\mathbf{X}}\boldsymbol{\eta})V_1) \boldsymbol{\eta}(U_1)]\\ \nonumber
&& [(\pounds_{\mathbf{X}}g)(h'U_1,V_1)+g((\pounds_{\mathbf{X}}h')U_1,V_1)] =0.
\end{eqnarray}
Putting $V_1=\varsigma$ in (\ref{3.19}) and using (\ref{3.16}),(\ref{3.17}) gives
\begin{equation}
(\verb"k"+1) g(\pounds_\mathbf{X} \varsigma, U_1) = 0,
\end{equation}
which implies either $\verb"k"+1=0$ or, $\verb"k"+1 \neq 0$.\par

Case I: When $\verb"k"+1=0$: Dileo and Pastore\cite{dil1} shown that in an $(\verb"k",\boldsymbol{\mu})'$-akm if $\verb"k"+1=0$, then $h'=0$ and hence the manifold is locally a warped product $\mathbf{N}' \times_{f} \mathbf{M}^{2m}$.\par

Case II: If $\verb"k"+1 \neq 0$, hence $g(\pounds_\mathbf{X} \varsigma,U_1) =0$, then (\ref{3.16}) implies
\begin{equation}\nonumber
(\pounds_\mathbf{X} \boldsymbol{\eta})U_1 = 0,
\end{equation}
which means that $\mathbf{X}$ is a strict infinitesimal contact transformation.\par
Hence we have:

\begin{theo}
If $\mathbf{N}^{2m+1}$ admits a Killing vector field $\mathbf{X}$, then either $N^{2m+1}$ is locally a warped product $\mathbf{N}' \times_{f} \mathbf{M}^{2m}$ or $\mathbf{X}$ is a strict infinitesimal contact transformation.
\end{theo}

\vspace{.9cm}
 Since $\mathbf{X}$ is Killing, $\rho=0$. Also, from Proposition 2.1, we have $r=2m(\verb"k"-2m)$. Hence, (\ref{1.3}) and (\ref{1.4}) implies
\begin{equation}\label{a.1}
 (\pounds_{\mathbf{X}}R)(U_1,V_1)W_1=0\,\, {\rm and} \,\, (\pounds_{\mathbf{X}}S)(U_1,V_1)=0.
\end{equation}
Taking Lie-derivative of (\ref{2.6}) and using (\ref{a.1}) , we obtain
\begin{eqnarray}\label{a.2}
R(U_1,V_1)\pounds_{\mathbf{X}}{\varsigma} &=& \verb"k"[((\pounds_{\mathbf{X}}\boldsymbol{\eta})V_1)U_1-((\pounds_{\mathbf{X}}\boldsymbol{\eta})U_1)V_1]\\ \nonumber
&&-\boldsymbol{\mu}[((\pounds_{\mathbf{X}}\boldsymbol{\eta})V_1)h'U_1-((\pounds_{\mathbf{X}}\boldsymbol{\eta})U_1)h'V_1\\ \nonumber  &&+\boldsymbol{\eta}(V_1)(\pounds_{\mathbf{X}}h')U_1-\boldsymbol{\eta}(U_1)(\pounds_{\mathbf{X}}h')V_1].
\end{eqnarray}
Since, $\mathbf{X}$ is Killing, (\ref{1.1}) implies
\begin{equation}\label{a.3}
g(\nabla_{\varsigma}{\mathbf{X}},U_1)+g(\nabla_{U_1}{\mathbf{X}},\varsigma) =0,
\end{equation}
which entails
\begin{equation}\label{a.4}
(\pounds_{\mathbf{X}}\boldsymbol{\eta})U_1=g(\pounds_{\mathbf{X}}\varsigma,U_1),
\end{equation}
for any vector field $U_1$. Let $\{e_0=\varsigma,e_1,e_2,...,e_{2m}\}$ be an orthonormal basis of $\mathbf{N}^{2m+1}$. Then from the above equation, we have
\begin{equation}\label{a.5}
\sum_{i=0}^{2n}g(e_i,V_1)(\pounds_{\mathbf{X}}\boldsymbol{\eta})e_i =g(\pounds_{\mathbf{X}}\varsigma,V_1),\,\, \sum_{i=0}^{2n}g(e_i,h'V_1)(\pounds_{\mathbf{X}}\boldsymbol{\eta})e_i =g(\pounds_{\mathbf{X}}\varsigma,h'V_1) .
\end{equation}
Again, for a non-Kenmotsu $(\verb"k",\boldsymbol{\mu})'$-akm, it is known that
\begin{eqnarray}\label{a.6}
(\nabla_{U_1}h')V_1 &=&g((\verb"k"+1)U_1-h'U_1,V_1)\varsigma + \boldsymbol{\eta}(V_1)((\verb"k"+1)U_1-h'U_1)\\ \nonumber
&&-2(\verb"k"+1)\boldsymbol{\eta}(U_1)\boldsymbol{\eta}(V_1)\varsigma
\end{eqnarray}
for any $U_1,\,V_1$ and $\boldsymbol{\mu}=-2$.\par
Also, from \cite{wana}, we get
\begin{equation}\label{a.7}
tr(\pounds_{\mathbf{X}}h')=0.
\end{equation}
Contracting (\ref{a.2}) and using (\ref{a.4}), (\ref{a.5}) and (\ref{a.7}), we obtain
\begin{equation}\label{a.8}
S(V_1,\pounds_{\mathbf{X}}\varsigma)=2mk g(\pounds_{\mathbf{X}}\varsigma,V_1)+2g(\pounds_{\mathbf{X}}\varsigma,h'V_1)+2g(\pounds_{\mathbf{X}}h'V_1,\varsigma).
\end{equation}
From (\ref{a.3}), we have $g(\pounds_{\mathbf{X}}h'V_1,\varsigma)=-g(\pounds_{\mathbf{X}}\varsigma,h'V_1)$. Hence the foregoing equation implies
\begin{equation}\label{a.9}
S(V_1,\pounds_{\mathbf{X}}\varsigma)=2mk g(\pounds_{\mathbf{X}}\varsigma,V_1).
\end{equation}
Replacing $U_1$ by $\pounds_{\mathbf{X}}\varsigma$ in (\ref{3.17.1}) and using (\ref{a.4}) and (\ref{a.9}), we get
\begin{equation}\label{a.10}
(\verb"k"+1)\pounds_{\mathbf{X}}\varsigma +h'\pounds_{\mathbf{X}}\varsigma =0.
\end{equation}
Operating $h'$ in (\ref{a.10}) and using (\ref{2.7}), we provide
\begin{equation}\label{a.11}
h'\pounds_{\mathbf{X}}\varsigma -\pounds_{\mathbf{X}}\varsigma=0.
\end{equation}
The above two equations implies $(\verb"k"+2)\pounds_{\mathbf{X}}\varsigma=0$, which implies either $\verb"k"+2=0$ or $\verb"k"+2\neq0$.\par
Case I: If $\verb"k"=-2$, then it is a non-Kenmotsu $(\verb"k",-2)'$-akm and is locally isometric to the Riemannian product $\mathbb{H}^{m+1}(-4)\times \mathbb{R}^m$ \cite{dil1}.\par
Case II:  If $\verb"k"\neq-2$, then $\pounds_{\mathbf{X}}\varsigma=0$, hence from (\ref{a.4}), we get $(\pounds_{\mathbf{X}}\boldsymbol{\eta})U_1=0$.
Hence we have:
\begin{theo}
If a non-Kenmotsu $(\verb"k",\boldsymbol{\mu})'$-akm $\mathbf{N}^{2m+1}$ permits a Killing vector field $\mathbf{X}$, it is either $\mathbf{N}^{2m+1}$ locally isometric to the product space $\mathbb{H}^{m+1}(-4)\times \mathbb{R}^m$ or $\mathbf{X}$ is a strict infinitesimal contact transformation.

\end{theo}

\vspace{.9cm}

\section{$\boldsymbol{\eta}$-Ricci-Yamabe solitons with conformal vector fields}
The concept of $\boldsymbol{\eta}$-Ricci-Yamabe soliton comes as a generalization of Ricci-Yamabe soliton and is defined by \cite{sidd}
\begin{equation}\label{1.7}
	\pounds_{\mathbf{X}} g + 2\alpha_1 S + (2\lambda_1 - \beta_1 r)g + 2\nu_1 \boldsymbol{\eta}\otimes\boldsymbol{\eta} = 0,
\end{equation}
where $\alpha_1, \beta_1, \nu_1$ are constant. If $\nu_1 =0$, then $\boldsymbol{\eta}$-Ricci-Yamabe soliton reduces to a Ricci-Yamabe soliton and for $\nu_1\neq 0$, it is called proper $\boldsymbol{\eta}$-Ricci-Yamabe soliton.
This soliton turns into \par
(i) $\boldsymbol{\eta}$-Ricci soliton if $\alpha_1 =1, \beta_1=0$,\par
(ii) $\boldsymbol{\eta}$-Yamabe soliton if $\alpha_1 =0, \beta_1 = 1$,\par
(iii) $\boldsymbol{\eta}$-Einstein soliton if $\alpha_1 =1, \beta_1 = -1$.\par
Several authors have studied $\boldsymbol{\eta}$-Ricci solitons and Ricci-Yamabe solitons, including (\cite{bla5}, \cite{bla6}, \cite{desde}, \cite{10}, \cite{103}) and many others.\par

Assume that the $(\verb"k",\boldsymbol{\mu})'$-akm $\mathbf{N}^{2m+1}$ admits an $\boldsymbol{\eta}$-Ricci-Yamabe soliton with $CVF$. 
 Then from (\ref{1.7}) we have
\begin{equation}\label{5.1}
(\pounds_{\mathbf{X}}g)(U_1,V_1)+2\alpha_1 S(U_1,V_1)+(2\lambda_1-\beta_1 r)g(U_1,V_1)+2\nu_1 \boldsymbol{\eta}(U_1)\boldsymbol{\eta}(V_1)=0.
\end{equation}
Using (\ref{1.1}) in (\ref{5.1}), we get
\begin{equation}\label{5.2}
\rho g(U_1,V_1) + \alpha_1 S(U_1,V_1)+(\lambda_1-\frac{\beta_1}{2} r)g(U_1,V_1)+\nu_1 \boldsymbol{\eta}(U_1)\boldsymbol{\eta}(V_1)=0.
\end{equation}
Setting $U_1=V_1=\varsigma$ in (\ref{5.2}) entails that
\begin{equation}\label{5.3}
\rho=-2mk\alpha_1 -\lambda_1-\nu_1+\frac{\beta_1}{2}r,
\end{equation}
which is a constant. Again, contracting $U_1$ and $V_1$ in (\ref{5.2}), we infer
\begin{equation}\label{5.4}
\rho(2m+1)=-\lambda_1(2m+1)+\frac{\beta_1}{2}r(2m+1)-\nu_1-\alpha_1 r.
\end{equation}
In view of (\ref{5.3}) and (\ref{5.4}), we obtain
\begin{equation}\label{5.5}
\lambda_1 = \beta_1 m(\verb"k"-2m)-\rho+2m \alpha_1
\end{equation}
and
\begin{equation}\label{5.6}
\nu_1 = -2m\alpha_1(\verb"k"+1).
\end{equation}
It is known that in a $(\verb"k",\boldsymbol{\mu})'$-akm, 
$r=2m(\verb"k"-2m)$, a constant. Hence with the help of (\ref{1.5}) and (\ref{5.3}), we get
$\rho r=0$, which implies $\rho =0$.  Hence, equation (\ref{5.2}) implies
\begin{equation}
\alpha_1 S(U_1,V_1) = (\frac{\beta_1}{2}r-\lambda_1)g(U_1,V_1) -\nu_1 \boldsymbol{\eta}(U_1)\boldsymbol{\eta}(V_1),
\end{equation}
which represents an $\boldsymbol{\eta}$-Einstein manifold. Since, the coefficients are constant, from Proposition 2.4, it becomes an Einstein manifold. Therefore we have:
\begin{theo}
If a $(\verb"k",\boldsymbol{\mu})'$-akm admits an $\boldsymbol{\eta}$-Ricci-Yamabe soliton with $CVF$, then it becomes an Einstein manifold and the constants $\lambda_1$ and $\nu_1$ are given by
\begin{equation} \nonumber
\lambda_1 = \beta_1 m(\verb"k"-2m)-\rho+2m \alpha_1
\end{equation}
and
\begin{equation}\nonumber
\nu_1 = -2m\alpha_1(\verb"k"+1).
\end{equation}
\end{theo}

In particular, if we take $\alpha_1 =1$, $\beta_1 =0$. Then we have
\begin{equation} \nonumber
\lambda_1 = 2m-\rho\,\,\, {\rm and}\,\,\, \nu_1 = -2m(\verb"k"+1).
\end{equation}
Hence we have:
\begin{cor}
If a $(\verb"k",\boldsymbol{\mu})'$-akm permits an $\boldsymbol{\eta}$-Ricci soliton with $CVF$, then the constants $\lambda_1$ and $\nu_1$ are given by
\begin{equation} \nonumber
\lambda_1 = 2m-\rho\,\,\, {\rm and}\,\,\, \nu_1 = -2m(\verb"k"+1).
\end{equation}
\end{cor}

In particular, if we take $\alpha_1 =0$, $\beta_1 =1$. Then we have
\begin{equation} \nonumber
\lambda_1 = m(\verb"k"-2m)-\rho\,\,\, {\rm and}\,\,\, \nu_1 = 0.
\end{equation}
Hence we have:
\begin{cor}
	If a $(\verb"k",\boldsymbol{\mu})'$-akm concedes an $\boldsymbol{\eta}$-Yamabe soliton with conformal vector field, then the soliton becomes Yamabe soliton and $\lambda_1 = 2m-\rho$.
\end{cor}

In particular, if we take $\alpha_1 =1$, $\beta_1 =-1$. Then we have
\begin{equation} \nonumber
\lambda_1 = -m(\verb"k"-2m)-\rho+2m\,\,\, {\rm and}\,\,\, \nu_1 = -2m(\verb"k"+1).
\end{equation}
Hence we have:
\begin{cor}
	If a $(\verb"k",\boldsymbol{\mu})'$-akm admits an $\boldsymbol{\eta}$-Einstein soliton with $CVF$, then the constants $\lambda_1$ and $\nu_1$ are given by
\begin{equation} \nonumber
\lambda_1 = -m(\verb"k"-2m)-\rho+2m\,\,\, {\rm and}\,\,\, \nu_1 = -2m(\verb"k"+1).
\end{equation}
\end{cor}

Also, if we take $\nu_1=0$, then the equation (\ref{5.6}) implies $\verb"k"=-1$ hence $h=0$. Therefore, the manifold becomes a Kenmotsu manifold. Hence we have:
\begin{cor}
If a $(\verb"k",\boldsymbol{\mu})'$-akm permits a proper Ricci-Yamabe soliton with $CVF$, then the manifold becomes a Kenmotsu manifold.
\end{cor}

\section{Example}

We choose $\mathbf{N}^{3} = \lbrace (x,y,z)\in \mathbb{R}^3\rbrace$, in which $(x,y,z)$ denote the standard coordinates in $\mathbb{R}^3$. Let $\varsigma, \varepsilon _2, \varepsilon _3$ be three vector fields in $\mathbb{R}^3$ which satisfies \cite{dil1}
\begin{equation}\nonumber
[\varsigma,\varepsilon _2] = -\varepsilon _2 - \varepsilon _3,\hspace{.3cm} [\varsigma,\varepsilon _3] = -\varepsilon _2 -\varepsilon _3,\hspace{.3cm} [\varepsilon _2,\varepsilon _3] = 0.
\end{equation}

The Riemannian metric $g$ is defined by
\begin{equation}\nonumber
g(\varsigma,\varsigma) = g(\varepsilon _2,\varepsilon _2) = g(\varepsilon _3,\varepsilon _3) = 1 \hspace{.3cm} and \hspace{.3cm} g(\varsigma,\varepsilon _2) = g(\varsigma,\varepsilon _3) = g(\varepsilon _2,\varepsilon _3) = 0.
\end{equation}

Here, $\boldsymbol{\eta}$ and $\varphi$ are defined as follows:\par
 $\boldsymbol{\eta}(W_1) = g(W_1,\varsigma)$, for any $W_1 \in \chi(\mathbf{N})$,
and
\begin{equation}\nonumber
\varphi \varsigma = 0,\hspace{.3cm} \varphi \varepsilon _2 = \varepsilon _3,\hspace{.3cm} \varphi \varepsilon _3 = -\varepsilon _2.
\end{equation}
Making use of the linearity property of $\varphi$ and $g$, we acquire
\begin{equation}\nonumber
\boldsymbol{\eta}(\varsigma) = 1,
\end{equation}
\begin{equation}\nonumber
\varphi^2 U_1 = -U_1 + \boldsymbol{\eta}(U_1)\varsigma,
\end{equation}
\begin{equation}\nonumber
g(\varphi U_1, \varphi V_1) = g(U_1,V_1) - \boldsymbol{\eta}(U_1)\boldsymbol{\eta}(V_1)
\end{equation}
for any $U_1, V_1 \in \chi(\mathbf{N})$. Thus the structure $(\varphi,\varsigma,\boldsymbol{\eta},g)$ is an almost contact structure.

Moreover, $h'\varsigma = 0,\hspace{.3cm} h'\varepsilon _2 = \varepsilon _3 \hspace{.3cm} and\hspace{.3cm} h'\varepsilon _3 = \varepsilon _2.$

In \cite{man}, the authors acquired the expression of the curvature tensor as:

\begin{equation}\nonumber
R(\varsigma,\varepsilon _2)\varsigma = 2(\varepsilon _2 + \varepsilon _3), \hspace{.3cm} R(\varsigma,\varepsilon _2)\varepsilon _2 = -2\varsigma,\hspace{.3cm} R(\varsigma,\varepsilon _2)\varepsilon _3 = -2\varsigma ,
\end{equation}
\begin{equation}\nonumber
R(\varepsilon _2,\varepsilon _3)\varsigma = R(\varepsilon _2,\varepsilon _3)\varepsilon _2 = R(\varepsilon _2,\varepsilon _3)\varepsilon _3 = 0,
\end{equation}
\begin{equation}\nonumber
R(\varsigma,\varepsilon _3)\varsigma = 2(\varepsilon _2 + \varepsilon _3),\hspace{.3cm} R(\varsigma,\varepsilon _3)\varepsilon _2 = -2\varsigma,\hspace{.3cm} R(\varsigma,\varepsilon _3)\varepsilon _3 = -2\varsigma.
\end{equation}
From the expressions of the $R$, we ensure that with $\verb"k"=-2$ and $\boldsymbol{\mu} = -2$, 
$\varsigma$ belongs to the $(\verb"k",\boldsymbol{\mu})'$-nullity distribution.\par
Foregoing expressions reveal the values of the Ricci tensor as:
\begin{equation}\nonumber
S(\varsigma,\varsigma) = -4,\hspace{.3cm}  S(\varepsilon _2,\varepsilon _2) = S(\varepsilon _3,\varepsilon _3) = -2.
\end{equation}
Hence, 
$r = S(\varsigma,\varsigma)+S(\varepsilon _2,\varepsilon _2) = S(\varepsilon _3,\varepsilon _3) =-8$.
Therefore, from (\ref{5.2}) we obtain
\begin{eqnarray}\nonumber
&&\rho-4\alpha_1+\lambda_1+4\beta_1+\nu_1=0\\ \nonumber
{\rm and} && \rho -2\alpha_1+\lambda_1+4\beta_1=0.
\end{eqnarray}
From the above two equations, we get
\begin{equation}\nonumber
\lambda_1 = 2\alpha_1-4\beta_1-\rho\,\, {\rm and}\,\, \nu_1 =2\alpha_1.
\end{equation}
Therefore, the manifold defines an $\boldsymbol{\eta}$-Ricci-Yamabe soliton for $\lambda_1 = 2\alpha_1-4\beta_1-\rho\,\, {\rm and}\,\, \nu_1 =2\alpha_1$.

\end{document}